\documentclass{amsart}
\usepackage[latin1]{inputenc}
\usepackage{graphicx}
\usepackage{latexsym}
\usepackage{amssymb}
\usepackage{amscd}

\newtheorem{thm}{Theorem}
\newtheorem{cor}[thm]{Corollary}
\newtheorem{lem}[thm]{Lemma}
\newtheorem{prop}[thm]{Proposition}
\theoremstyle{definition}
\newtheorem{defn}[thm]{Definition}
\theoremstyle{remark}

\newtheorem{prob}{Problem}
\numberwithin{equation}{section}
\newcommand{\To}{\longrightarrow}

\begin{document}
\setcounter{tocdepth}{1}
%\date{\today}
\title[]{A Boolean algebra and a Banach space obtained by push-out iteration}
\author{Antonio Avil\'es and Christina Brech}
\address{Universidad de Murcia, Departamento de Matem\'{a}ticas, Campus de Espinardo 30100 Murcia, Spain} \email{avileslo@um.es}
\thanks{A. Avil\'{e}s was supported by MEC and FEDER (Project MTM2008-05396), Fundaci\'{o}n S\'{e}neca
(Project 08848/PI/08), Ram\'{o}n y Cajal contract (RYC-2008-02051) and an FP7-PEOPLE-ERG-2008 action.}
\address{Departamento de Matem\'{a}tica, Instituto de Matem\'{a}tica e Estat\'{\i}stica, Universidade de S\~{a}o Paulo, Rua do Mat\~ao, 1010 - 05508-090, S\~ao Paulo, Brazil}
\thanks{C. Brech was supported by FAPESP grant (2007/08213-2), which is part of Thematic Project FAPESP (2006/02378-7).}
\email{christina.brech@gmail.com}

\subjclass{Primary 06E05; Secondary 03E35, 03G05, 46B26, 54G05}

\begin{abstract}
Under the assumption that $\mathfrak c$ is a regular cardinal, we prove the existence and uniqueness of a Boolean algebra $\mathfrak B$ of size $\mathfrak c$ defined by sharing the main structural properties that $\mathcal{P}(\omega)/fin$ has under CH and in the $\aleph_2$-Cohen model. We prove a similar result in the category of Banach spaces.
\end{abstract}

\maketitle

\section{Introduction}

In this paper two lines of research converge, one related to the theory of Boolean algebras and the other one to Banach spaces.\\

In the context of Boolean algebras, the topic goes back to Parovi\v{c}enko's theorem~\cite{Parovicenko}, which establishes that, under CH, $\mathcal{P}(\omega)/fin$ is the unique Boolean algebra of size $\mathfrak c$ with the property that given any diagram of embeddings of Boolean algebras like
$$\begin{CD}
 S\\
 @AAA \\
 R @>>>\mathcal{P}(\omega)/fin,
\end{CD}$$
where $R$ and $S$ are countable, there exists an embedding $S\To\mathcal{P}(\omega)/fin$ which makes the diagram commutative. This characterization is indeed equivalent to CH~\cite{vanDouwenvanMill}. There has been a line of research~\cite{DowHart,FucGesSheSou,FucGesSou,FucKopShe,Geschke,Steprans} showing that many results about $\mathcal{P}(\omega)/fin$ under CH can be generalized to the $\aleph_2$-Cohen model\footnote{By ($\kappa$-)Cohen model, we mean a model obtained by adding ($\kappa$ many) Cohen reals to a model of CH} (and to a less extent to any Cohen model). The key point
was proven by Stepr\={a}ns~\cite{Steprans} and is that in that case $\mathcal{P}(\omega)/fin$ is tightly $\sigma$-filtered. Later, Dow and Hart~\cite{DowHart} introduced the notion of Cohen-Parovi\v{c}enko Boolean algebras. On the one hand, $\mathcal{P}(\omega)/fin$ is Cohen-Parovi\v{c}enko in Cohen models, and on the other hand there exists a unique Cohen-Parovi\v{c}enko Boolean algebra of size $\mathfrak c$ whenever $\mathfrak c\leq\aleph_2$. Therefore this property characterizes $\mathcal{P}(\omega)/fin$ in the $\aleph_2$-Cohen model in the spirit of Parovi\v{c}enko's theorem.\\

Our main result states that whenever $\mathfrak c$ is a regular cardinal, there exists a unique tightly $\sigma$-filtered Cohen-Parovi\v{c}enko Boolean algebra.  Essentially, we are proving that whenever $\mathfrak c$ is regular, there exists a canonical Boolean algebra that behaves like $\mathcal{P}(\omega)/fin$ does under CH and in the $\aleph_2$-Cohen model.
By a result of Geschke~\cite{Geschke}, this Boolean algebra cannot be $\mathcal{P}(\omega)/fin$ when $\mathfrak c>\aleph_2$. Hence, in the $\kappa$-Cohen model ($\kappa >\aleph_2$), $\mathcal{P}(\omega)/fin$ and our Boolean algebra are two non-isomorphic Cohen-Parovi\v{c}enko Boolean algebras, and this solves a question of Dow and Hart~\cite[Questions 2 and 3]{DowHart} about uniqueness for this property.\\

We found that many notions being used in the literature on this topic can be reformulated using the concept of push-out from category theory. Apart from aesthetic considerations, this approach has an objective advantage: it allows us to export ideas to other categories like we do with Banach spaces, where we have push-outs with similar properties.\\

In the context of Banach spaces, an analogue of Parovi\v{c}enko's theorem has been established by Kubi\'{s}~\cite{Kubis}: Under CH, there exists a unique Banach space $\mathfrak X$ of density $\mathfrak c$ with the property that given any diagram of isometric embeddings like
$$\begin{CD}
 S\\
 @AAA \\
 R @>>>\mathfrak X
\end{CD}$$
where $R$ and $S$ are separable, there exists an embedding $S\To\mathfrak X$ which makes the diagram commutative. The existence of such a space can be proven in ZFC but not its uniqueness~\cite{extremadura}. In this paper we shall consider a stronger property than the one stated above which will imply existence and uniqueness of the space $\mathfrak X$ under the weaker assumption that $\mathfrak c$ is a regular cardinal.\\

In both situations the method of construction is the same as in \cite{DowHart} and \cite{extremadura} respectively, by a long chain of push-outs. The difficulty lies in proving that the algebra or the space constructed in this way is indeed unique. It would be nice to have a unified approach for Banach spaces and Boolean algebras in the context of category theory, but our attempts to do so seemed to become too technical and to obscure both subjects rather than to enlighten them. The reader will find, nevertheless, an obvious parallelism.\\

The paper is structured as follows: In Section~\ref{sectionBoole} we prove the existence and uniqueness of our Boolean algebra $\mathfrak B$. In Section~\ref{sectionBanach} we do the same with our Banach space $\mathfrak X$. The two sections run parallel but they can be read independently, and for the convenience of the reader they include a self-contained account of the facts that we need about push-outs in each category. In Section~\ref{sectionStone} we consider the compact space $\mathfrak K$, Stone dual of $\mathfrak B$, we establish that $\mathfrak X$ is isometric to a subspace of the space of continuous functions $C(\mathfrak K)$ and we also prove that $\mathfrak K$ is homogeneous with respect to $P$-points, generalizing the results of Rudin~\cite{Rudin} under CH, Stepr\={a}ns~\cite{Steprans} in the $\aleph_2$-Cohen model, and Geschke~\cite{Geschke}. In Section~\ref{sectionProblems} we state several open problems and in Section~\ref{sectionCardinals} we point out more general versions of our results for different cardinals.\\

We would like to thank Wies\l aw Kubi\'{s} for some relevant remarks that helped to improve this paper.

\section{Boolean algebras}\label{sectionBoole}

\subsection{Preliminary definitions} The push-out is a general notion of category theory, and the one that we shall use here refers to the category of Boolean algebras. We shall consider only push-outs made of embeddings (one-to-one morphisms), although it is a more general concept. For this reason, we present the subject in a different way than usual, more convenient for us and equivalent for the case of embeddings.\\

The join and meet operations in a Boolean algebra $B$ are denoted by $\vee$ and $\wedge$, and the complement of $r\in B$ by $\overline{r}$. The subalgebra generated by $H$ is $\langle H\rangle$.

\begin{defn}
Let $B$ be a Boolean algebra and let $S$ and $A$ be subalgebras of $B$. We say that $B$ is the internal push-out of $S$ and $A$ if the following conditions hold:
\begin{enumerate}
\item $B=\langle S\cup A\rangle$.
\item For every $a\in A$ and every $s\in S$, if $a\wedge s = 0$, then there exists $r\in A\cap S$ such that $a\leq r$ and $s\leq \overline{r}$.\\
\end{enumerate} 
\end{defn}

We notice that condition (2) above can be substituted by the following equivalent one:
\begin{itemize}
\item[(2')] For every $a\in A$ and every $s\in S$, if $a\leq s$, then there exists $r\in A\cap S$ such that $a\leq r\leq s$.\\
\end{itemize}

The same definition can be found in \cite{Geschke} as ``$A$ and $S$ commute''. 

Suppose that we have a diagram of embeddings of Boolean algebras:

$$\begin{CD}
 S @>>> B\\
 @AAA @AAA\\
 R @>>>A.
\end{CD}$$

We say that it is a push-out diagram if, when all algebras are viewed as subalgebras of $B$, we have that $R=S\cap A$ and $B$ is the internal push-out of $A$ and $S$.\\

If $B$ is the push-out of $S$ and $A$, then $B$ is isomorphic to $(S\otimes A)/V$ where $S\otimes A$ denotes the free sum of $S$ and $A$, and $V$ is the ideal generated by the formal intersections $r\wedge \overline{r}$, where $r\in A\cap S$ (viewing $r\in S$, $\overline{r}\in A$). Also, given a diagram of embeddings

$$\begin{CD}
 S \\
 @AAA\\
 R @>>>A,
\end{CD}$$
there is a unique way (up to isomorphism) to complete it into a push out diagram, putting $B = (S\otimes A)/V$ where $V$ is as above.\\

If $B$ is the push-out of $S$ and $A$ we will write $B = \mathbf{PO}[S,A]$. Sometimes we want to make explicit the intersection space $R=S\cap A$, and then we write $B=\mathbf{PO}_R[S,A]$, meaning that $B$ is the push-out of $S$ and $A$ and $R=S\cap A$.\\

\begin{defn} An embedding of Boolean algebras $A\To B$ is said to be a posex (push-out separable extension) if there exists a push-out diagram of embeddings
$$\begin{CD}
 A @>>> B\\
 @AAA @AAA\\
 R @>>>S
\end{CD}$$
such that $S$ and $R$ are countable.
\end{defn}

 If $A\subset B$, we can rephrase this definition saying that $B$ is a posex of $A$ if there exists a countable subalgebra $S\subset B$ such that $B=\mathbf{PO}[S,A]$.\\

\begin{defn}
We say that a Boolean algebra $\mathfrak B$ is tightly $\sigma$-filtered if there exists an ordinal $\lambda$ and a family $\{B_\alpha : \alpha\leq\lambda\}$ of subalgebras of $\mathfrak B$ such that
\begin{enumerate}
\item $B_\alpha\subset B_\beta$ whenever $\alpha<\beta\leq\lambda$,
\item $B_0 = \{0,1\}$ and $B_\lambda = \mathfrak B$,
\item $B_{\alpha+1}$ is a posex of $B_\alpha$ for every $\alpha<\lambda$,
\item $B_\beta = \bigcup_{\alpha<\beta}B_\alpha$ for every limit ordinal $\beta\leq\lambda$.\\
\end{enumerate}
\end{defn}

We shall see later that the definition of Koppelberg~\cite{Koppelberg} of tightly $\sigma$-filtered Boolean algebra is equivalent to this one.

\subsection{The main result}

We are now ready to state the main result of this section.

\begin{thm}[$\mathfrak c$ is a regular cardinal]\label{existsuniqueB}
There exists a unique (up to isomorphism) Boolean algebra $\mathfrak B$ with the following properties:
\begin{enumerate}
\item $|\mathfrak B| = \mathfrak c$,
\item $\mathfrak B$ is tightly $\sigma$-filtered,
\item For any diagram of embeddings of the form
$$\begin{CD}
 B \\
 @AAA \\
 A @>>> \mathfrak B,
\end{CD}$$
if $|A|<\mathfrak c$ and $A\To B$ is posex, then there exists an embedding $B\To \mathfrak B$ which makes the diagram commutative.
\end{enumerate}
\end{thm} 

\subsection{Basic properties of push-outs}

Before going to the proof of Theorem~\ref{existsuniqueB}, we collect some elementary properties of push-outs of Boolean algebras.

\begin{prop}\label{propertiesPOBoole} The following facts about push-outs hold:
\begin{enumerate}
\item If $B=\mathbf{PO}[S_\alpha,A]$ for every $\alpha$, where $\{S_\alpha : \alpha<\xi\}$ is an increasing chain of Boolean subalgebras of $B$, $S_0\subset S_1\subset\cdots$, then $B=\mathbf{PO}[\bigcup_\alpha S_\alpha,A]$.

\item If $B=\mathbf{PO}[S,A]$ and we have $S'\subset S$ and $A'\subset A$, then $B = \mathbf{PO}[\langle S\cup A'\rangle, \langle A\cup S'\rangle]$.

\item If $B_1 = \mathbf{PO}_{S_0}[S_1,B_0]$ and $B_2 = \mathbf{PO}_{S_1}[S_2,B_1]$, then $B_2 = \mathbf{PO}_{S_0}[S_2,B_0]$. 

\item Suppose that we have two increasing sequences of subalgebras of $B$, $S_0\subset S_1\subset\cdots$ and $B_0\subset B_1\subset\cdots$, and $B_{n} = \mathbf{PO}_{S_0}[S_n,B_0]$ for every $n$. Then  $\bigcup_n B_n = \mathbf{PO}_{S_0}[\bigcup_n S_n, B_0]$.

\item Suppose that we have two increasing sequences of subalgebras of $B$, $S_0\subset S_1\subset\cdots$ and $B_0\subset B_1\subset\cdots$, and $B_{n+1} = \mathbf{PO}_{S_{n}}[S_{n+1},B_n]$ for every $n$. Then  $\bigcup_n B_n = \mathbf{PO}_{S_0}[\bigcup_n S_n, B_0]$.

\item Suppose that $A=\bigcup_{i\in I}A_i\subset B$, $S = \bigcup_{j\in J}S_j \subset B$ and $\langle A_i\cup S_j\rangle = \mathbf{PO}[S_j,A_i]$ for every $i$ and every $j$. Then $\langle A \cup S\rangle = \mathbf{PO}[S,A]$.

\end{enumerate}
\end{prop}

Proof: (1) is trivial.

For item (2) it is enough to consider the case where $S'=\emptyset$, that is to prove $B = \mathbf{PO}[\langle S\cup A'\rangle,A]$, the general case being obtained by the successive application of the cases $S'=\emptyset$ and its symmetric one $A'=\emptyset$. Clearly $B=\langle S\cup A \rangle = \langle (S\cup A')\cup A \rangle$. Let $s \in S \cup A'$ and $a \in A$ be such that $a \leq s$. If $s \in S$, our hypothesis guarantees that there is $r \in A \cap S \subseteq A \cap (S \cup A')$ such that $a \leq r \leq s$. If $s \in A' \subseteq A \cap (S \cup A')$, take $r = s$ and we are done.

For item (3), it is clear that $S_2 \cap B_0  = S_0$ and $B_2 = \langle S_2\cup B_0\rangle$. On the other hand assume that $a\leq s_2$ for some $a\in B_0$ and $s\in S_2$. Then, since $a\in B_1$ and  $B_2 = \mathbf{PO}_{S_1}[S_2,B_1]$, there exists $r_1\in S_1$ such that $a\leq r_1\leq s$. Since we also have $B_1 = \mathbf{PO}_{S_0}[S_1,B_0]$, we find $r_0\in S_0$ with $a\leq r_0\leq r_1\leq s$.

Item (4) is evident (it holds even for transfinite sequences though we do not need to use that) and item (5) follows from combining (3) and (4). Item (6) is also easy to see.$\qed$\\

We remark that some of the properties of push-outs of Banach spaces in Proposition~\ref{propertiesPO} do not hold in the context of Boolean algebras. Namely, suppose that $A$ is freely generated by $\{x_i : i\in I\}$ and $B$ is freely generated by $\{x_i : i\in I\}\cup\{y\}$. If $S=\langle y\rangle$, then $B=\mathbf{PO}[S,A]$. However, if we consider $D=\langle x_i, x_i \wedge y : i\in I\rangle$, then $A\subset D\subset B$ but $B\neq\mathbf{PO}[S,D]$ and $D\neq \mathbf{PO}[S\cap D,A]$. So, we don't have an analogue of items (1) and  (2) of Proposition~\ref{propertiesPO}.

\begin{prop}\label{characterizationposexBoole}
$B$ is a posex of $A$ if and only if the following two conditions hold:
\begin{enumerate}
\item $B$ is countably generated over $A$ (that is, $B = \langle A\cup Q\rangle$ for some countable set $Q$), and
\item For every $b\in B$, the set $\{a\in A : a\leq b\}$ is a countably generated ideal of~$A$.
\end{enumerate}
\end{prop}

Proof: We suppose first that $B$ is a posex of $A$, $B = \mathbf{PO}_R[S,A]$ for some countable $S$. It is obvious that $B$ is countably generated over $A$ so we concentrate in proving the other property. It is clear that property (2) holds whenever $b\in A$ (in that case the ideal is just generated by $b$) and also when $b\in S$, since by the push-out property the ideal  $I_b = \{a\in A : a\leq b\}$ is generated by $I_b\cap R$, which is countable. Let us consider now an arbitrary element $b\in B$. We can express it in the form $b = (s_1 \wedge a_1) \vee\cdots\vee (s_n \wedge a_n)$ where $a_1,\ldots,a_n$ form a partition in $A$ (that is, they are disjoint and their join is 1) and $s_1,\ldots,s_n\in S$. Pick now $a\in A$ such that $a\leq b$. Then $a \wedge a_i \leq s_i$ for every $i$, so since $a \wedge a_i\in A$ and $s_i\in S$, there exists $r_i\in R$ with $a \wedge a_i \leq r_i \leq s_i$. We have that
$$ a = (a \wedge a_1) \vee\cdots\vee (a \wedge a_n) \leq (r_1 \wedge a_1) \vee\cdots\vee (r_n \wedge a_n) \leq (s_1 \wedge a_1) \vee\cdots\vee (s_n \wedge a_n) = b .$$ 
Notice that $(r_1 \wedge a_1) \vee\cdots\vee (r_n \wedge a_n) \in A$. It follows that the ideal $\{a\in A : a\leq b\}$ is generated by the countable set
$$\{ r = (r_1 \wedge a_1) \vee\cdots\vee (r_n \wedge a_n) : r_i\in R, r\leq b\}$$ 
We prove now the converse implication. So we assume that $B=\langle A\cup Q\rangle$ with $Q$ countable, and for every $b\in B$, we fix a countable set $G_b\subset A$ that generates the ideal $\{a\in A : a\leq b\}$. We define an increasing sequence of subalgebras of $B$ making $S_0 = \langle Q\rangle$ and $S_{n+1} = \langle S_n\cup\bigcup_{b\in S_n}G_b\rangle$. All these are countably generated -hence countable- Boolean algebras. Taking $S= \bigcup_{n<\omega}S_n$ we get that $B = \mathbf{PO}[S,A]$. Namely, if $a\leq s$ with $s\in S$ and $a\in A$, then $s\in S_n$ for some $n$, and then there exists $r\in G_s$ with $a\leq r\leq s$. Just observe that $r\in G_s\subset A\cap S_{n+1}\subset A\cap S$.$\qed$\\

Condition (2) of Proposition~\ref{characterizationposexBoole} is found in the literature under the following names: $A$ is a good subalgebra of $B$ \cite{Steprans}; $A$ is an $\aleph_0$-ideal subalgebra of $B$~\cite{DowHart}; $A$ is a $\sigma$-subalgebra of $B$ \cite{FucKopShe,Geschke}. We keep the latter terminology. So, $B$ is a posex of $A$ if and only if $A$ is a $\sigma$-subalgebra of $B$ and $B$ is countably generated over $A$. Conversely, $A$ is a $\sigma$-subalgebra of $B$ if and only if every intermediate algebra $A\subset C \subset B$ which is countably generated over $A$ is a posex of $A$. After these equivalences, it becomes obvious that Koppelberg's definition~\cite{Koppelberg} of tightly $\sigma$-filtered algebra and our own are the same, cf. also \cite[Theorem 2.5]{Geschke}.\\

\begin{prop}\label{propertiesposexB} The following facts hold:
\begin{enumerate}
\item If $A\subset B$ and $B$ is countable, then $B$ is a posex of $A$.
\item If $B$ is a posex of $A$, $A\subset B'\subset B$, and $B'$ is countably generated over $A$, then $B'$ is a posex of $A$.
\item If $B_1$ is a posex of $B_0$ and $B_2$ is a posex of $B_1$, then $B_2$ is a posex of $B_0$.
\item If we have $B_n\subset B_{n+1}$ and $B_n$ is a posex of $B_0$ for every $n<\omega$, then $\bigcup_n B_n$ is a posex of $B_0$.
\item If $B_{n+1}$ is a posex of $B_n$ for every $n\in\mathbb{N}$, then $\bigcup_n B_n$ is a posex of $B_0$.
\item If $B$ is a posex of $A$ and $S_0\subset B$ is countable, then there exists a countable subalgebra $S_0\subset S\subset B$ with $B = \mathbf{PO}[S,A]$.
\item If we have $B_0\subset B_1\subset B_2$, $B_2$ is a posex of $B_0$ and $B_1$ is countably generated over $B_0$, then $B_2$ is a posex of $B_1$.

\end{enumerate}
\end{prop}

Proof: (1) is trivial, (2) follows from Proposition~\ref{characterizationposexBoole}. We prove now (3). Suppose that $B_1 = \mathbf{PO}_{R_0}[S_0,B_0]$ and $B_2 = \mathbf{PO}_{T_0}[U_0,B_1]$. Our objective is to apply Proposition~\ref{propertiesPOBoole}(3), so we need to overcome the difficulty that $S_0\neq T_0$. Inductively on $n$, we will define countable subalgebras $R_n, S_n, T_n, U_n$ forming four increasing sequences so that  $B_1 = \mathbf{PO}_{R_n}[S_n,B_0]$ and $B_2 = \mathbf{PO}_{T_n}[U_n,B_1]$ for every $n$. The inductive procedure is as follows: pick a countable set $Q_n\subset B_0$ such that $T_n\subset\langle S_n \cup Q_n\rangle$, then define $S_{n+1} = \langle S_n\cup Q_n\rangle$, $R_{n+1} = S_{n+1}\cap B_0$, $U_{n+1} = \langle U_n\cup S_{n+1}\rangle$ and $T_{n+1} = U_{n+1}\cap B_1$. The push-out relations $B_1 = \mathbf{PO}_{R_n}[S_n,B_0]$ and $B_2 = \mathbf{PO}_{T_n}[U_n,B_1]$ follow from  Proposition~\ref{propertiesPOBoole}(2). Also notice that $T_n\subset S_{n+1}\subset T_{n+1}$ for every $n$. Hence $\bigcup_n S_n = \bigcup_n T_n$. On the other hand, by Proposition~\ref{propertiesPOBoole}(4) we have that   $B_1 = \mathbf{PO}_{\bigcup_n R_n}[\bigcup_n S_n,B_0]$ and $B_2 = \mathbf{PO}_{\bigcup_n T_n}[\bigcup_n U_n,B_1]$ . By Proposition~\ref{propertiesPOBoole}(3),  $B_2 = \mathbf{PO}_{\bigcup_n R_n}[\bigcup_n U_n,B_0]$, which proves that $B_2$ is a posex of $B_0$, since $\bigcup_n U_n$ is countable.

Item (4) is proven easily using Proposition~\ref{characterizationposexBoole}, and (5) is a consequence of (3) and (4). For (6), consider first a countable subalgebra $S_1\subset B$ such that $B = \mathbf{PO}[S_1,A]$. Then find a countable set $Q\subset A$ such that $S_0\subset \langle Q\cup S_1\rangle$. From Proposition~\ref{propertiesPOBoole}(2) we get that $B = \mathbf{PO}[\langle Q\cup S_1\rangle,A]$. To prove (7), use (2) to get that $B_1$ is a posex of $B_0$ hence $B_1 = \mathbf{PO}[S_1,B_0]$ for some countable subalgebra $S_1$. By (6) there exists a countable subalgebra $S_2$ such that $B_2 = \mathbf{PO}[S_2,B_0]$ and $S_1\subset S_2$. By Proposition~\ref{propertiesPOBoole}(2), $B_2 = \mathbf{PO}[S_2,\langle B_0\cup S_1\rangle]$ but $\langle B_0\cup S_1\rangle = B_1$.$\qed$\\

Consider again the example where $A$ is freely generated by $\{x_i : i\in I\}$ with $I$ uncountable, $B$ is freely generated by $\{x_i : i\in I\}\cup\{y\}$ and $D=\langle x_i, x_i \wedge y : i\in I\rangle$. Then $B=\mathbf{PO}[\langle y\rangle,A]$, so $B$ is a posex of $A$. On the other hand, $A\subset D\subset B$ but neither $D$ is a posex of $A$ nor $B$ is a posex of $D$. Indeed $D$ is not countably generated over $A$, and $\{a\in D : a\leq y\}$ is not a countably generated ideal of $D$.\\

\subsection{Proof of Theorem~\ref{existsuniqueB}}

The following definition, as well as Lemma~\ref{POidealB} are due to Geschke~\cite{Geschke}. A proof of Lemma~\ref{POidealB} can also be obtained by imitating the proof of Lemma~\ref{POidealE}.

\begin{defn}
Let $\mathfrak B$ be a (uncountable) Boolean algebra. An additive $\sigma$-skeleton of $\mathfrak B$ is a family $\mathcal{F}$ of subalgebras of $\mathfrak B$ with the following properties:
\begin{enumerate}
\item $\{0,1\}\in\mathcal{F}$
\item For every subfamily $\mathcal{G}\subset\mathcal{F}$, we have $\left\langle \bigcup\mathcal{G}\right\rangle\in\mathcal{F}$.
\item For every infinite subalgebra $A\subset \mathfrak B$, there exists $A'\in\mathcal{F}$ such that with $A\subset A'\subset \mathfrak B$ and $|A| = |A'|$.
\item Every $A\in \mathcal{F}$ is a $\sigma$-subalgebra of $\mathfrak B$.\\ 
\end{enumerate}
\end{defn}

We will often use the following property of an additive $\sigma$-skeleton: If $A\in \mathcal F$, $A\subset A'\subset\mathfrak B$ and $A'$ is countably generated over $A$, then there exists $B\in\mathcal F$ such that $A'\subset B$ and $B$ is countably generated over $A'$. This is a direct consequence of properties (2) and (3): Suppose that $A' = \langle A\cup S\rangle$ with $S$ countable Boolean algebra. Then there exists a countable $S_1\in \mathcal F$ with $S\subset S_1$, and we can take $B = \langle A\cup S_1\rangle$.

\begin{lem}[Geschke]\label{POidealB}
For a Boolean algebra $\mathfrak B$ the following are equivalent:
\begin{enumerate}
\item $\mathfrak B$ is tightly $\sigma$-filtered.
\item There exists an additive $\sigma$-skeleton $\mathcal{F}$ of $\mathfrak B$.
\end{enumerate}
\end{lem}

We can now prove Theorem~\ref{existsuniqueB}. First we prove existence. We consider $\mathfrak c = \bigcup_{\alpha<\frak c} \Phi_\alpha$ a decomposition of the continuum into $\mathfrak c$ many subsets of cardinality $\mathfrak c$ such that $\alpha\leq \min(\Phi_\alpha)$ for every $\alpha$. We define recursively an increasing chain of Boolean algebras $\{B_\alpha : \alpha<\mathfrak c\}$, so that at the end $\mathfrak B=\bigcup_{\alpha<\mathfrak c} B_\alpha$. We start with $B_0 = \{0,1\}$. After $B_\alpha$ is defined, we consider a family $\{(R_\gamma,S_\gamma) : \gamma\in \Phi_\alpha\}$ such that
\begin{itemize}
\item For every $\gamma \in \Phi_\alpha$, $S_\gamma$ is a countable Boolean algebra and $R_\gamma = B_\alpha\cap S_\gamma$.
\item For every countable subalgebra $R\subset B_\alpha$ and every countable superalgebra $S\supset R$ there exists $\gamma\in\Phi_\alpha$ and a Boolean isomorphism $j:S\To S_\gamma$ such that $R=R_\gamma$ and $j(x) = x$ for $x\in R$.
\end{itemize}

For limit ordinals $\beta$ we define $B_\beta = \bigcup_{\alpha<\beta}B_\alpha$. At successor stages we define $B_{\alpha+1} = \mathbf{PO}_{R_\alpha}[S_\alpha,B_\alpha]$. By construction, it is clear that $\mathfrak B$ is a tightly $\sigma$-filtered Boolean algebra of cardinality $\mathfrak c$. We check property (3) in the statement of the theorem. Suppose that we have $A\subset \mathfrak B$ with $|A|<\mathfrak c$, and $B=\mathbf{PO}_R[S,A]$ with $S$ countable. By the regularity of $\mathfrak c$, we can find $\alpha<\mathfrak c$ such that $A\subset B_\alpha$. Then, there exists $\gamma\in \Phi_\alpha$ such that $R=R_\gamma$ and (modulo an isomorphism) $S=S_\gamma$, so that $B_{\gamma+1} =\mathbf{PO}_R[S,B_\gamma]$. Consider $\tilde{B} = \langle S\cup A\rangle \subset B_{\gamma+1}$, so that $\tilde{B} = \mathbf{PO}_R[S,A]$. Since push-out is unique up to isomorphism, we can find an isomorphism $\tilde{u}:B\To\tilde{B}\subset \mathfrak B$ such that $\tilde{u}(a) =a$ for all $a\in A$.\\

We prove now uniqueness. Suppose that we have two Boolean algebras like this, $\mathfrak B$ and $\mathfrak B'$. We consider their respective additive $\sigma$-skeletons $\mathcal{F}$ and $\mathcal{F}'$ that witness that they are tightly $\sigma$-filtered. Let us suppose that $\mathfrak B=\{x_\alpha : \alpha<\mathfrak c\}$ and $\mathfrak B' = \{y_\alpha : \alpha<\mathfrak c\}$. We shall construct recursively two increasing chains of subalgebras $\{B_\alpha : \alpha<\mathfrak c\}$ and $\{B'_\alpha : \alpha<\mathfrak c\}$ and a family of Boolean isomorphisms $f_\alpha:B_\alpha \To B'_\alpha$ with the following properties:
\begin{enumerate}
\item The isomorphisms are coherent, that is $f_\beta|_{B_\alpha} = f_\alpha$ whenever $\alpha<\beta$.
\item For every $\alpha$, $B_\alpha\in\mathcal{F}$ and $B'_\alpha\in\mathcal{F}'$.
\item $x_\alpha\in B_{\alpha+1}$ and $y_\alpha\in B'_{\alpha+1}$. In this way we make sure that $\mathfrak B=\bigcup_{\alpha<\mathfrak c} B_\alpha$ and $\mathfrak B' = \bigcup_{\alpha<\mathfrak c} B'_\alpha$.
\item $B_{\alpha+1}$ is countably generated over $B_\alpha$ for every $\alpha$. This implies that $|B_\alpha| = |\alpha|$ for every $\alpha\geq \omega$.
\end{enumerate}
After this, the isomorphisms $f_\alpha$ induce a global isomorphism $f:\mathfrak B\To \mathfrak B'$. We proceed to the inductive construction. We start with $B_0 =\{0,1\}$ and $B'_0=\{0,1\}$. If $\beta$ is a limit ordinal, we simply put $B_\beta = \bigcup_{\alpha<\beta}B_\alpha$, $B'_\beta = \bigcup_{\alpha<\beta}B'_\alpha$ an the isomorphism $f_\beta$ is induced by the previous ones. Now we show how to construct $B_{\alpha+1}$, $B'_{\alpha+1}$ and $f_{\alpha+1}$ from the previous ones. We construct inductively on $n$, sequences of subalgebras $B_\alpha[n]\subset \mathfrak B$ , $B'_\alpha[n]\subset \mathfrak B'$ and coherent isomorphisms $f_\alpha[n]:B_\alpha[n]\To B'_\alpha[n] $ as in the picture:

$$\begin{array}[c]{cccccccccccc}
 B_\alpha & = & B_\alpha[0] & \subset & B_{\alpha}[1] &\subset & B_{\alpha}[2] & \subset & B_{\alpha}[3] & \cdots &\subset & \mathfrak B \\
          &   &  \downarrow &         & \downarrow    &        & \downarrow    &         & \downarrow    &        &        & \\
 B'_\alpha & = & B'_\alpha[0] & \subset & B'_{\alpha}[1] &\subset & B'_{\alpha}[2] & \subset & B'_{\alpha}[3] & \cdots &\subset & \mathfrak B' \\
\end{array}$$
and we will make $B_{\alpha+1} = \bigcup_{n<\omega}B_\alpha[n]$, $B'_{\alpha+1} = \bigcup_{n<\omega}B'_\alpha[n]$ and $f_{\alpha+1}$ induced by the isomorphisms $f_\alpha[n]$. All the algebras $B_\alpha[n+1]$ and $B'_\alpha[n+1]$ will be countably generated over $B_\alpha[n]$ and $B'_\alpha[n]$ respectively. The inductive procedure is as follows. There are two cases:\\

Case 1: $n$ is even. Then, we define $B_\alpha[n+1]$ to be such that $x_\alpha\in B_\alpha[n+1]$, $B_\alpha[n+1]$ is countably generated over $B_\alpha[n]$,  and $B_\alpha[n+1]\in \mathcal{F}$. Since $B_\alpha\in \mathcal{F}$ which is an additive $\sigma$-skeleton, $B_\alpha[n+1]$ is a posex of $B_\alpha$, and by Proposition~\ref{propertiesposexB}(7), also $B_\alpha[n+1]$ is a posex of $B_\alpha[n]$. Hence, since $\mathcal B'$ satisfies the statement of our theorem, we can find a Boolean embedding $f_\alpha[n+1]: B_\alpha[n+1]\To \mathcal B'$ such that $f_\alpha[n+1]|_{B_\alpha[n]} = f_\alpha[n]$. We define finally $B'_\alpha[n+1] = f_\alpha[n+1](B_\alpha[n+1])$.\\

Case 2: $n$ is odd. Then, we define $B'_\alpha[n+1]$ to be such that $y_\alpha\in B'_\alpha[n+1]$, $B'_\alpha[n+1]$ is countably generated over $B'_\alpha[n]$  and $B'_\alpha[n+1]\in \mathcal{F}'$. Since $B'_\alpha\in \mathcal{F}'$ which is an additive $\sigma$-skeleton, $B'_\alpha[n+1]$ is a posex of $B_\alpha$, hence also a posex of $B'_\alpha[n]$, so since $\mathfrak B$ satisfies the statement of our theorem, we can find an embedding $g_\alpha[n+1]: B'_\alpha[n+1]\To \mathfrak B$ such that $g_\alpha[n+1]|_{B_\alpha[n]} = f^{-1}_\alpha[n]$. We define finally $B_\alpha[n+1] = g_\alpha[n+1](B_\alpha[n+1])$ and $f_\alpha[n+1] = g^{-1}_\alpha[n+1]$.\\

Proceeding this way, we have that $B_\alpha[n]\in\mathcal{F}$ for $n$ odd, while $B'_\alpha[n]\in \mathcal{F}'$ for $n$ even. At the end, $B_{\alpha+1} = \bigcup_{n<\omega}B_\alpha[2n+1]\in\mathcal{F}$ and $B'_{\alpha+1} = \bigcup_{n<\omega}B'_\alpha[2n]\in\mathcal{F}'$, which concludes the proof. $\qed$\\

\subsection{Remarks} Dow and Hart~\cite{DowHart} say that $B$ is $(*,\aleph_0)$-ideal, if for every $\kappa<\mathfrak c$ there exists $\kappa$-cub of $\sigma$-subalgebras of $B$. If we have an additive $\sigma$-skeleton $\mathcal{F}$ for $B$, then the algebras in $\mathcal{F}$ of cardinality $\kappa$ form a $\kappa$-cub of $\sigma$-subalgebras of $B$. Hence, every tightly $\sigma$-filtered Boolean algebra of cardinality $\mathfrak c$ is $(*,\aleph_0)$-ideal.\\

Also, they say that a subalgebra $A\subset B$ is $\aleph_0$-ideal complete if for every two orthogonal countably generated ideals $I,J$ of $A$, there exists $c\in B$ such that $I = \{a\in A: a\leq c\}$ and $J = \{a\in A: a\leq \bar{c}\}$. We have the following fact:

\begin{prop}
$A$ is an $\aleph_0$-ideal complete subalgebra of $B$ if and only if for every posex $A\To C$ there exists an embedding $g:C\To B$ with $g|_A = 1_A$.
\end{prop}

Proof: If the statement in the right holds, and we take two orthogonal countably generated ideals $I,J\subset A$, we can consider $R$ a countable subalgebra of $A$ generated by the union of countable sets of generators of $I$ and $J$. Take a superalgebra of the form $S=\langle R\cup\{c\}\rangle$ where $x\leq c$ for every $x\in I$ and $x\leq \bar{c}$ for every $x\in J$, and let $C = \mathbf{PO}_R[S,A]$. Our assumption provides an embedding $g:C\To B$ and the element $g(c)$ is the desired one.

Assume now that $A$ is $\aleph_0$-ideal complete and consider $f:A\To C$ posex. It is enough to consider the case when $C$ is finitely generated over $f(A)$, so suppose that $C= \langle f(A)\cup\{c_1,c_2,\ldots,c_n\} \rangle$ where $\{c_1,\ldots,c_n\}$ form a partition. For every $i<n$ find $d_i\in B$ such that $$\{a\in A : f(a)\leq c_i\} = \{a\in A: a\leq d_i\}$$ $$\{a\in A : f(a)\leq \overline{c_i}\} = \{a\in A: a\leq \overline{d_i}\}$$
Define $d_n=1$. Convert the $d_i$'s into a partition by setting $d'_k = d_k \setminus (\bigvee_{i<n}d_i)$ for $k\leq n$. Then it is possible to define the desired extension $g$ by declaring $g(c_k) = d'_k$.$\qed$\\

Finally, Dow and Hart call a Boolean algebra $\mathfrak B$ to be Cohen-Parovi\v{c}enko if:
\begin{enumerate}
\item $|\mathfrak{B}|=\mathfrak c$.
\item $\mathfrak B$ is $(\ast,\aleph_0)$-ideal.
\item All subalgebras of $\mathfrak B$ of cardinality less than $\mathfrak c$ are $\aleph_0$-ideal complete.
\end{enumerate}

Dow and Hart prove that, if $\mathfrak c\leq\aleph_2$, there is a unique Cohen-Parovi\v{c}enko Boolean algebra of size $\mathfrak c$ and they ask if this is true in ZFC. The answer is negative. On the one hand, $\mathcal{P}(\omega)/fin$ is Cohen-Parovi\v{c}enko in any Cohen model, as they show based on results of Stepr\={a}ns~\cite{Steprans}. On the other hand, as it follows from the comments above, items (1) and (3) are equivalent to the corresponding items of Theorem~\ref{existsuniqueB}, while (2) is weaker, which guarantees that the Boolean algebra of Theorem~\ref{existsuniqueB} is Cohen-Parovi\v{c}enko. 
As a consequence of a result of Geschke we get the following result, which guarantees that the Boolean algebra of Theorem~\ref{existsuniqueB} is not $\mathcal{P}(\omega)/fin$, whenever $\mathfrak c > \aleph_2$:

\begin{thm}
If $\mathfrak B = \mathcal{P}(\omega)/fin$ is the algebra of Theorem~\ref{existsuniqueB}, then $\mathfrak c \leq \aleph_2$.
\end{thm}

Proof: Geschke~\cite{Geschke} proves that a complete Boolean algebra of size greater than $\aleph_2$ cannot be tightly $\sigma$-filtered. Hence, if $\mathfrak c >\aleph_2$, $\mathcal{P}(\omega)$ is not tightly $\sigma$-filtered. This implies that neither $\mathcal{P}(\omega)/fin$ is such, because a tower of subalgebras witnessing tight $\sigma$-filtration could be lifted to $\mathcal{P}(\omega)$.$\qed$\\ 

Another property of the algebra $\mathfrak B$ of Theorem~\ref{existsuniqueB} is:

\begin{prop}\label{Buniversal}
If $A$ is a tightly $\sigma$-filtered Boolean algebra with $|A|\leq \frak c$, then $A$ is isomorphic to a subalgebra of $\mathfrak B$.
\end{prop}

Proof: Let $\{A_\alpha: \alpha<\kappa\}$ be subalgebras that witness that $A$ is tightly $\sigma$-filtered. We can suppose that $\kappa$ is the cardinality of $A$, cf.\cite{Geschke}. Then, inductively we can extend a given embedding $A_\alpha\To \mathfrak B$ to $A_{\alpha+1}\To \mathfrak{B}$, by the properties of $\mathfrak B$.$\qed$\\

\section{Banach spaces}\label{sectionBanach}

\subsection{Preliminary definitions} The push-out is a general notion of category theory, and the one that we shall use here refers to the category $Ban_1$ of Banach spaces, together with operators of norm at most 1. We shall consider only push-outs made of isometric embeddings of Banach spaces, although it is a more general concept. For this reason, we present the subject in a different way than usual, more convenient for us but equivalent for the case of isometric embeddings.  

\begin{defn}
Let $Y$ be a Banach space and let $S$ and $X$ be subspaces of $Y$. We say that $Y$ is the internal push-out of $S$ and $X$ if the following conditions hold:
\begin{enumerate}
\item $Y=\overline{S+X}$
\item $\|x+s\| = \inf\{\|x+r\|+\|s-r\| : r\in S\cap X\}$ for every $x\in X$ and every $s\in S$.\\
\end{enumerate} 
\end{defn}

Suppose that we have a diagram of isometric embeddings of Banach spaces like:

$$\begin{CD}
 S @>>> Y\\
 @AAA @AAA\\
 R @>>>X.
\end{CD}$$

We say that it is a push-out diagram if, when all spaces are seen as subspaces of $Y$, we have that $R=S\cap X$ and $Y$ is the internal push-out of $X$ and $S$.\\

Note that $Y$ being the push-out of $S$ and $X$ means that $Y$ is isometric to the quotient space $(S\oplus_{\ell_1} X)/V$ where $V = \{(r,-r) : r\in X\cap S\}$. In particular, $Y=S+X = \{s+x : s\in S, x\in X\}$. Also, given a diagram of isometric embeddings

$$\begin{CD}
 S \\
 @AuAA\\
 R @>v>>X,
\end{CD}$$
there is a unique way (up to isometries) to complete it into a push out diagram, precisely by making $Y = (S\oplus_{\ell_1} X)/\tilde{R}$ and putting the obvious arrows, where $\tilde{R} = \{(u(r),v(-r)) : r\in R\}$.\\

If $Y$ is the push-out of $S$ and $X$ as above, we will write $Y = \mathbf{PO}[S,X]$. Sometimes we want to make explicit the intersection space $R=S\cap X$, and then we write $Y=\mathbf{PO}_R[S,X]$.

\begin{defn} An isometric embedding of Banach spaces $X\To Y$ is said to be a posex if there exists a push-out diagram of isometric embeddings
$$\begin{CD}
 X @>>> Y\\
 @AAA @AAA\\
 R @>>>S
\end{CD}$$
such that $S$ and $R$ are separable.
\end{defn}

 If $X\subset Y$, we can rephrase this definition saying that there exists a separable subspace $S\subset Y$ such that $Y=\mathbf{PO}[S,X]$.\\

\begin{defn}
We say that a Banach space $\mathfrak X$ is tightly $\sigma$-filtered if there exists an ordinal $\lambda$ and a family $\{X_\alpha : \alpha\leq\lambda\}$ of subspaces of $\mathfrak X$ such that
\begin{enumerate}
\item $X_\alpha\subset X_\beta$ whenever $\alpha<\beta\leq\lambda$,
\item $X_0 = 0$ and $X_\lambda = X$,
\item $X_{\alpha+1}$ is a posex of $X_\alpha$ for every $\alpha<\lambda$,
\item $X_\beta = \overline{\bigcup_{\alpha<\beta}X_\alpha}$ for every limit ordinal $\beta\leq\lambda$.\\
\end{enumerate}
\end{defn}

\subsection{The main result}

\begin{thm}[$\mathfrak c$ is a regular cardinal]\label{existsuniqueE}
There exists a unique (up to isometry) Banach space $\mathfrak X$ with the following properties:
\begin{enumerate}
\item $dens(\mathfrak X) = \mathfrak c$,
\item $\mathfrak X$ is tightly $\sigma$-filtered,
\item For any diagram of isometric embeddings of the form
$$\begin{CD}
 Y \\
 @AAA \\
 X @>>> \mathfrak X,
\end{CD}$$
if $dens(X)<\mathfrak c$ and $X\To Y$ is posex, then there exists an isometric embedding $Y\To \mathfrak X$ which makes the diagram commutative.
\end{enumerate}
\end{thm} 

\subsection{Basic properties of push-outs}

Before entering the proof of Theorem~\ref{existsuniqueE}, we collect some elementary properties of push-outs that we shall use. For the sake of completeness, we include their proofs.

\begin{prop}\label{propertiesPO} The following facts about push-outs hold:
\begin{enumerate}
\item Suppose $Y = \mathbf{PO}[S,X]$, $S\subset S'\subset Y$ and $X\subset X'\subset Y$. Then $Y = \mathbf{PO}[S',X']$.
\item Suppose $Y = \mathbf{PO}[S,X]$, and $X\subset Y'\subset Y$. Then $Y' = \mathbf{PO}[S\cap Y',X]$.
\item Suppose $Y = \mathbf{PO}[S,X]$, and $S\subset Y'\subset Y$. Then $Y' = \mathbf{PO}[S,X\cap Y']$.
\item Let $X\subset Y\subset Z$ and $S_0\subset S_1\subset S_2\subset Z$ be such that $Y = \mathbf{PO}_{S_0}[S_1,X]$ and $Z = \mathbf{PO}_{S_1}[S_2,Y]$. Then $Z = \mathbf{PO}_{S_0}[S_2,X]$.
\item Let $X_0\subset X_1\subset\cdots$ and $S_0\subset S_1\subset\cdots$ be Banach spaces such that for every $n\in\mathbb{N}$ we have that $X_{n+1} = \mathbf{PO}_{S_n}[S_{n+1},X_n]$. Let $X=\overline{\bigcup_n X_n}$ and $S=\overline{\bigcup_n S_n}$. Then $X=\mathbf{PO}_{S_0}[S,X_0]$.
\item\label{amalgamation} Suppose $X = \overline{\bigcup_{i\in I}X_i}\subset Y$, $S = \overline{\bigcup_{j\in J}S_j}\subset Y$ and $\overline{X_i+S_j} = \mathbf{PO}[S_j,X_i]$ for every $i$ and every $j$. Then $\overline{X+S} = \mathbf{PO}[S,X]$.
\end{enumerate}
\end{prop}

Proof: For (1) it is enough to consider the case when $X=X'$. The case $S=S'$ is just the same, and the general fact follows from the application of those two. So we have to prove that given $x\in X$, $s'\in S'$ and $\varepsilon>0$ there exists $t\in X\cap S'$ such that
$$\|x+s'\| +\varepsilon > \|x+t\| + \|s'-t\|,$$
Write $s'=s+x'$ where $s\in S$ and $x'\in X$. Then there exists $r\in X \cap S$ such that
$$\|x+s'\| = \|x+x'+s\|>\|x+x'+r\|+\|s-r\|-\varepsilon = \|x+(x'+r)\| + \|s'-(x'+r)\|-\varepsilon$$
so just take $t=x'+r$.

Item (2) is immediate, simply notice that $Y' = X+(Y'\cap S)$, because if $y\in Y'$ and $y=x+s$ with $x\in X$ and $s\in S$, then $s=y-x\in Y'\cap S$. Item (3) is symmetric to item (2).

For (4), consider $\varepsilon>0$, $x\in X$ and $s_2\in S_2$. Then since $Z=\mathbf{PO}_{S_1}[S_2,Y]$, there exists $s_1\in S_1$ such that
$$ \|x+s_2\| > \|x+s_1\|+\|s_2-s_1\|-\varepsilon/2$$  
and since $Y=\mathbf{PO}_{S_0}[S_1,X]$ there exists $s_0\in S_0$ with

$$ \|x+s_1\| > \|x+s_0\|+\|s_1-s_0\|-\varepsilon/2$$
so finally

$$ \|x+s_2\| > \|x+s_0\|+\|s_1-s_0\|+\|s_2-s_1\|-\varepsilon \geq \|x+s_0\|+\|s_2-s_0\|-\varepsilon$$ 

For (5), we know by repeated application of item (4) that $X_{n+1} = \mathbf{PO}_{S_0}[S_n,X_0]$. Consider $\varepsilon>0$, $x_0\in X_0$ and $s\in S$. Then there exists $n$ and $s_n\in S_n$ such that $\|s_n-s\|<\varepsilon/3$. Then, we can find $r\in S_0$ such that
$$\|x_0+s\| > \|x_0+s_n\|-\varepsilon/3  > \|x_0+r\| + \|s_n-r\|-2\varepsilon/3 >  \|x_0+r\| + \|s-r\|-\varepsilon$$
Item (6) is proven similarly by approximation.$\qed$\\

\begin{cor}\label{propertiesposex} The following facts about posexes hold:
\begin{enumerate}
\item If $Y$ is a posex of $X$, then $Y/X$ is separable.
\item $X\subset Y$ and $Y$ is separable, then $Y$ is a posex of $X$.
\item If $Y$ is a posex of $X$ and $X\subset Y'\subset Y$, then $Y'$ is a posex of $X$.
\item If $Y$ is a posex of $X$ and $X\subset X'\subset Y$, then $Y$ is a posex of $X'$.
\item If $Y$ is a posex of $X$ and $Z$ is a posex of $Y$, then $Z$ is a posex of $X$.
\item If $X_{n+1}$ is a posex of $X_n$ for every $n\in\mathbb{N}$, then $\overline{\bigcup_n X_n}$ is a posex of $X_0$. 
\item\label{sequenceposex} If we have $X_n\subset X_{n+1}$ and $X_n$ is a posex of $X_0$ for every $n\in\mathbb{N}$, then $\overline{\bigcup_n X_n}$ is a posex of $X_0$.
\end{enumerate}
\end{cor}

A subspace $X$ of a Banach space $\mathfrak X$ will be called a $\sigma$-subspace of $\mathfrak X$, if for every $X\subset Y\subset \mathfrak X$, if $Y/X$ is separable, then $Y$ is a posex of $X$.  

\begin{defn}
Let $\mathfrak{X}$ be a (nonseparable) Banach space. An additive $\sigma$-skeleton of $\mathfrak X$ is a family $\mathcal{F}$ of subspaces of $\mathfrak X$ with the following properties:
\begin{enumerate}
\item $0\in\mathcal{F}$
\item For every subfamily $\mathcal{G}\subset\mathcal{F}$, we have $\overline{span}\left(\bigcup\mathcal{G}\right)\in\mathcal{F}$.
\item For every infinite-dimensional $X\subset \mathfrak X$ there exists $Y\in\mathcal{F}$ with $X\subset Y\subset \mathfrak X$ and $dens(Y) = dens(X)$.
\item Each $X\in\mathcal{F}$ is a $\sigma$-subspace of $\mathfrak X$.\\ 
\end{enumerate}
\end{defn}

We will often use the following property of an additive $\sigma$-skeleton: If $X\in \mathcal F$, $X\subset Y\subset\mathfrak X$ and $Y/X$ is separable, then there exists $Z\in\mathcal F$ such that $Y\subset Z$ and $Z/Y$ is separable. The proof is a direct consequence of properties (2) and (3): Suppose that $Y = \overline{X+S}$ with $S$ separable. Then there exists a separable $S_1\in \mathcal F$ with $S\subset S_1$, and we can take $Z = \overline{X+S_1}$.

We prove now a result for Banach spaces corresponding to Lemma \ref{POidealB}.

\begin{lem}\label{POidealE}
For a Banach space $\mathfrak X$ the following are equivalent:
\begin{enumerate}
\item $\mathfrak X$ is tightly $\sigma$-filtered.
\item $\mathfrak X$ has an additive $\sigma$-skeleton.
\end{enumerate}
\end{lem}

Proof: That $(2)$ implies $(1)$ is evident: it is enough to define the subspaces $X_\alpha$ inductively just taking care that $X_\alpha\in\mathcal{F}$ for every $\alpha$. Now, we suppose that we have a tower of subspaces $\{X_\alpha : \alpha\leq\lambda\}$ like in $(2)$. For every $\alpha<\lambda$ we consider separable subspaces $R_\alpha\subset S_\alpha$ such that $X_{\alpha+1} = \mathbf{PO}_{R_\alpha}[S_\alpha,X_\alpha]$.\\

Given a set of ordinals $\Gamma\subset\lambda$, we define $E(\Gamma) = \overline{span}\left(\bigcup_{\gamma\in\Gamma}S_\gamma\right)$.
We say that the set $\Gamma$ is saturated if for every $\alpha\in\Gamma$ we have that $R_\alpha\subset E(\Gamma\cap\alpha)$. We shall prove that the family $\mathcal{F} = \{E(\Gamma) : \Gamma\subset \lambda \text{ is saturated}\}$ is the additive $\sigma$-skeleton that we are looking for.\\

It is clear that $0=E(\emptyset)$. Also, if we have a family $\{\Gamma_i : i\in I\}$ of saturated sets, then $\overline{span}\left(\bigcup_i E(\Gamma_i)\right) = E\left(\bigcup_i \Gamma_i\right)$, so the union of saturated sets is saturated. 

Given any countable set $\Gamma \subseteq \lambda$, it is possible to find a countable saturated set $\Delta$ such that $\Gamma \subset \Delta$: We define $\Delta = \{\delta_s: s \in \omega^{<\omega}\}$ where $(\delta_s)_{s\in \omega^{<\omega}}$ is defined inductively on the length of $s\in\omega^{<\omega}$ as follows. Let $\Gamma = \{\delta_{(n)} : n\in \omega\}$ and given $s\in \omega^{<\omega}$, let $\{\delta_{s^\frown n} : n\in \omega\}$ be such that $\delta_{s^\frown n}<\delta_s$ and $R_{\delta_s}\subset E(\{\delta_{s^\frown n} : n\in \omega\})$. Notice that $\Delta$ is a countable and saturated set which contains $\Gamma$ and also $\sup (\delta+1: \delta \in\Gamma) = \sup (\delta+1: \delta \in\Delta)$.\\

Now, suppose $X \subseteq \mathfrak X$ and let us find a saturated set $\Delta$ such that $X \subseteq E(\Delta)$ and $dens(X) = dens(E(\Delta))$. Since the union of saturated sets is saturated, we can assume without loss of generality that $X$ is separable and find a countable set $\Gamma$ such that $X\subset E(\Gamma)$. Take $\Delta$ as in the previous paragraph and notice that
$X \subset E(\Delta) \in \mathcal F$ and $E(\Delta)$ is separable since $\Delta$ is countable.\\

It remains to prove that if $X\in\mathcal{F}$ and we have $X\subset Y\subset \mathfrak X$ with $Y/X$ separable, then $Y$ is a posex of $X$. It is enough to prove the following statement, and we shall do it by induction on $\delta_1 = \sup(\delta+1 : \delta\in \Delta)$:\\

``For every saturated set $\Gamma\subset\lambda$ and every countable set $\Delta\subset\lambda$, $E(\Gamma\cup\Delta)$ is a posex of $E(\Gamma)$.''\\

We fix $\delta_1<\lambda$ and we assume that the statement above holds for all saturated sets $\Gamma\subset\lambda$ and for all countable sets $\Delta'\subset\lambda$ with $\sup(\delta+1 : \delta\in \Delta')<\delta_1$.\\

\emph{Case} 1: $\delta_1$ is a limit ordinal. This follows immediately from the inductive hypothesis 
using Corollary~\ref{propertiesposex}(\ref{sequenceposex}).\\

\emph{Case} 2: $\delta_1 = \delta_0+1$ for some $\delta_0$ and $\Delta = \{\delta_0\}$. We distinguish two subcases:\\

\emph{Case} 2a: $\sup(\Gamma)\leq\delta_0$. Consider the chain of subspaces

$$E(\Gamma) \subset \overline{E(\Gamma)+R_{\delta_0}}\subset E(\Gamma\cup\{\delta_0\}).$$

The left hand extension is a posex extension by the inductive hypothesis, because there exists a countable set $\Delta'$ such that $R_{\delta_0}\subset E(\Delta')$ and $\sup(\delta+1,\delta\in\Delta')\leq\delta_0<\delta_1$. The right hand extension is also posex, because we have the push-out diagram
$$\begin{CD}
 S_{\delta_0} @>>> X_{\delta_1}\\
 @AAA @AAA\\
 R_{\delta_0} @>>>X_{\delta_0}
\end{CD}$$
and since $\sup(\Gamma)\leq\delta_0$  (and $\delta_0\not\in\Gamma$, otherwise it is trivial), we can interpolate

$$\begin{CD}
 S_{\delta_0} @>>> E(\Gamma\cup\{\delta_0\}) @>>> X_{\delta_1}\\
 @AAA @AAA @AAA\\
 R_{\delta_0} @>>> \overline{E(\Gamma)+R_{\delta_0}} @>>> X_{\delta_0}
\end{CD}$$
where clearly $E(\Gamma\cup\{\delta_0\}) = \overline{E(\Gamma)+R_{\delta_0} + S_{\delta_0}}$, so that the left hand square is a push-out diagram.\\

\emph{Case} 2b: $\sup(\Gamma)>\delta_0$. For every $\xi\leq\lambda$, we call $\Gamma_\xi = \Gamma\cap\xi$. By Case 2a, there exists a separable space $S$ such that $$E(\Gamma_{\delta_1}\cup\{\delta_0\}) = \mathbf{PO}[S,E(\Gamma_{\delta_1})].$$

We can suppose that $S_{\delta_0}\subset S$. We shall prove by induction on $\xi$ that $$E(\Gamma_{\xi}\cup\{\delta_0\}) = \mathbf{PO}[S,E(\Gamma_{\xi})] \text{ for }\delta_1\leq\xi\leq\lambda.$$ 

If $\xi$ is a limit ordinal, then $E(\Gamma_\xi) = \overline{\bigcup_{\eta<\xi}E(\Gamma_\eta)}$, and we just need Proposition~\ref{propertiesPO}(\ref{amalgamation}). Otherwise suppose that $\xi = \eta+1$. In the nontrivial case, $\eta\in\Gamma$ and $\Gamma_\xi = \Gamma_{\eta}\cup\{\eta\}$. By the inductive hypothesis we have that $$E(\Gamma_{\eta}\cup\{\delta_0\}) = \mathbf{PO}[S,E(\Gamma_\eta)]$$ and on the other hand $X_\xi = X_{\eta+1} = \mathbf{PO}[S_\eta,X_\eta]$.

Suppose that we are given $x\in E(\Gamma_\xi)$, $s\in S$ and $\varepsilon>0$ and we have to find $r\in E(\Gamma_\xi)\cap S$ such that $\|x+s\|\geq \|x+r\| + \|s-r\| - \varepsilon$. We write $x= s_\eta + x_\eta$ where $s_\eta\in S_\eta$ and $x_\eta \in E(\Gamma_\eta)$. Notice that $x_\eta+s\in E(\Gamma_\eta)+S\subset E(\Gamma_\eta \cup \{\delta_0\})\subset X_\eta$ since $S_{\delta_0}\subset X_\eta$  ($\delta_0<\eta$ as $\delta_1<\xi$). Hence, using that $X_\xi = \mathbf{PO}[S_\eta,X_\eta]$ we find $r_\eta\in S_\eta\cap X_\eta$ such that 

$$\|s_\eta + x_\eta + s\| \geq \|x_\eta+s+r_\eta\| + \|s_\eta -r_\eta\| - \varepsilon/2.$$

Since $\eta\in\Gamma$ and $\Gamma$ is saturated, we have that $R_\eta\subset E(\Gamma_\eta)$, therefore $x_\eta+r_\eta\in E(\Gamma_\eta)$ and $s\in S$. Hence, using that $E(\Gamma_{\eta}\cup\{\delta_0\}) = \mathbf{PO}[S,E(\Gamma_\eta)]$ we get $r\in S\cap E(\Gamma_\eta)$ such that

$$\|x_\eta+s+r_\eta\| \geq \|s +r \| + \| x_\eta+r_\eta-r\| - \varepsilon/2.$$

Combining both inequalities,
$$ \|s_\eta + x_\eta + s\| \geq \|s +r \| + \| x_\eta+r_\eta-r\| + \|s_\eta -r_\eta\| - \varepsilon \geq \|s+r\| + \|x_\eta + s_\eta - r\| - \varepsilon,$$
as desired.\\

\emph{Case} 3: $\delta_1=\delta_0+1$ for some $\delta_0 \in \Delta$ and $|\Delta|>1$. We consider the set $\Delta\setminus\{\delta_0\}$. We have proven few paragraphs above in this proof that we can find a countable saturated set $\Delta'\supset\Delta\setminus\{\delta_0\}$ such that $\sup(\delta+1: \delta \in \Delta') \leq \delta_0$. By the inductive hypothesis $E(\Gamma\cup\Delta')$ is a posex of $E(\Gamma)$. Since $\Gamma \cup \Delta'$ is saturated, the already proven case 2 provides that $E(\Gamma\cup\Delta'\cup\{\delta_0\})$ is a posex of $E(\Gamma\cup\Delta')$. Composing both posex extensions we get that $E(\Gamma\cup\Delta'\cup\{\delta_0\})$ is a posex of $E(\Gamma)$. Since $E(\Gamma) \subset E(\Gamma\cup\Delta)\subset E(\Gamma\cup\Delta'\cup\{\delta_0\})$ we finally get that $E(\Gamma\cup\Delta)$ is a posex of $E(\Gamma)$.$\qed$\\

We finally prove Theorem~\ref{existsuniqueE}. First we prove existence. We consider $\mathfrak c = \bigcup_{\alpha<\frak c} \Phi_\alpha$ a decomposition of the continuum into $\mathfrak c$ many subsets of cardinality $\mathfrak c$ such that $\alpha\leq \min(\Phi_\alpha)$ for every $\alpha$. We define recursively an increasing chain of Banach spaces $\{X_\alpha : \alpha<\mathfrak c\}$, so that at the end $\mathfrak X=\bigcup_{\alpha<\mathfrak c} X_\alpha$. We start with $X_0 = 0$. After $X_\alpha$ is defined, we consider a family $\{(R_\gamma,S_\gamma) : \gamma\in \Phi_\alpha\}$ where
\begin{itemize}
\item For every $\gamma \in \Phi_\alpha$, $S_\gamma$ is a separable Banach space and $R_\gamma = X_\alpha\cap S_\gamma$.
\item For every separable subspace $R\subset X_\alpha$ and every separable superspace $S\supset R$ there exists $\gamma\in\Gamma$ and an isometry $j:S\To S_\gamma$ such that $R=R_\gamma$ and $j(x) = x$ for $x\in R$.
\end{itemize}

For limit ordinals $\beta$ we define $X_\beta = \overline{\bigcup_{\alpha<\beta}X_\alpha}$. At successor stages we define $X_{\alpha+1} = \mathbf{PO}_{R_\alpha}[S_\alpha,X_\alpha]$. By construction, it is clear that $\mathfrak X$ is a tightly $\sigma$-filtered Banach space of density $\mathfrak c$. We check property (3) in the statement of the theorem. Suppose that we have $X\subset \mathfrak X$ with $dens(X)<\mathfrak c$, and $Y=\mathbf{PO}_R[S,X]$ with $S$ separable. By the regularity of $\mathfrak c$, we can find $\alpha<\mathfrak c$ such that $X\subset X_\alpha$. Then, there exists $\gamma\in \Phi_\alpha$ such that $R=R_\gamma$ and (modulo an isometry) $S=S_\gamma$. Then $X_{\gamma+1} =\mathbf{PO}_R[S,X_\gamma]$. Consider $\tilde{Y} = S+X\subset X_{\gamma+1}$, so that $\tilde{Y} = \mathbf{PO}_R[S,X]$. Since push-out is unique up to isometry, we can find an isometry $\tilde{u}:Y\To\tilde{Y}\subset \mathfrak X$ such that $\tilde{u}(x) =x$ for all $x\in X$.\\

We prove now uniqueness. Suppose that we have two spaces like this, $\mathfrak X$ and $\mathfrak X'$. We consider their respective additive $\sigma$-skeletons $\mathcal{F}$ and $\mathcal{F}'$ that witness that they are tightly $\sigma$-filtered. Let us suppose that $\mathfrak X=\overline{span}\{x_\alpha : \alpha<\mathfrak c\}$ and $\mathfrak X' = \overline{span}\{y_\alpha : \alpha<\mathfrak c\}$. We shall construct recursively two increasing chains of subspaces $\{X_\alpha : \alpha<\mathfrak c\}$ and $\{X'_\alpha : \alpha<\mathfrak c\}$ and a family of bijective isometries $f_\alpha:X_\alpha \To X'_\alpha$ with the following properties:
\begin{enumerate}
\item The isometries are coherent, that is $f_\beta|_{X_\alpha} = f_\alpha$ whenever $\alpha<\beta$.
\item For every $\alpha$, $X_\alpha\in\mathcal{F}$ and $X'_\alpha\in\mathcal{F}'$.
\item $x_\alpha\in X_{\alpha+1}$ and $y_\alpha\in X'_{\alpha+1}$. In this way we make sure that $\mathfrak X=\bigcup_{\alpha<\mathfrak c}X_\alpha$ and $\mathfrak X' = \bigcup_{\alpha<\mathfrak c} X'_\alpha$.
\item Each quotient $X_{\alpha+1}/X_\alpha$ is separable. This implies that $dens(X_\alpha) = |\alpha|$ for every $\alpha\geq\omega$.
\end{enumerate}
After this, the isometries $f_\alpha$ induce a global isometry $f:\mathfrak X\To \mathfrak X'$. We proceed to the inductive construction. We start with $X_0 = 0$ and $X'_0=0$. If $\beta$ is a limit ordinal, we simply put $X_\beta = \overline{\bigcup_{\alpha<\beta}X_\alpha}$, $X'_\beta = \overline{\bigcup_{\alpha<\beta}X'_\alpha}$ and the isometry $f_\beta$ is induced by the previous isometries. Now, we show how to construct $X_{\alpha+1}$, $X'_{\alpha+1}$ and $f_{\alpha+1}$ from the previous ones. We construct inductively on $n$, sequences of subspaces $X_\alpha[n]\subset \mathfrak X$ , $X'_\alpha[n]\subset \mathfrak X'$ and coherent isometries $f_\alpha[n]:X_\alpha[n]\To X'_\alpha[n] $ as in the picture:

$$\begin{array}[c]{cccccccccccc}
 X_\alpha & = & X_\alpha[0] & \subset & X_{\alpha}[1] &\subset & X_{\alpha}[2] & \subset & X_{\alpha}[3] & \cdots &\subset & \mathfrak X \\
          &   &  \downarrow &         & \downarrow    &        & \downarrow    &         & \downarrow    &        &        & \\
 X'_\alpha & = & X'_\alpha[0] & \subset & X'_{\alpha}[1] &\subset & X'_{\alpha}[2] & \subset & X'_{\alpha}[3] & \cdots &\subset & \mathfrak X' \\
\end{array}$$
and we will make $X_{\alpha+1} = \overline{\bigcup_{n<\omega}X_\alpha[n]}$, $X'_{\alpha+1} = \overline{\bigcup_{n<\omega}X'_\alpha[n]}$ and $f_{\alpha+1}$ induced by the isometries $f_\alpha[n]$. All the quotient spaces $X_\alpha[n+1]/X_\alpha[n]$ and $X'_\alpha[n+1]/X'_\alpha[n]$ will be separable. The inductive procedure is as follows. There are two cases:\\

Case 1: $n$ is even. Then, we define $X_\alpha[n+1]$ to be such that $x_\alpha\in X_\alpha[n+1]$, $X_\alpha[n+1]/X_\alpha[n]$ is separable,  and $X_\alpha[n+1]\in \mathcal{F}$. Since $X_\alpha\in \mathcal{F}$ which is an additive $\sigma$-skeleton, $X_\alpha[n+1]$ is a posex of $X_\alpha[n]$, so since $\mathfrak X'$ satisfies the statement of our theorem, we can find an into isometry $f_\alpha[n+1]: X_\alpha[n+1]\To \mathfrak X'$ such that $f_\alpha[n+1]|_{X_\alpha[n]} = f_\alpha[n]$. We define finally $X'_\alpha[n+1] = f_\alpha[n+1](X_\alpha[n+1])$.\\

Case 2: $n$ is odd. Then, we define $X'_\alpha[n+1]$ to be such that $y_\alpha\in X'_\alpha[n+1]$, $X'_\alpha[n+1]/X'_\alpha[n]$ is separable,  and $X'_\alpha[n+1]\in \mathcal{F}'$. Since $X'_\alpha\in \mathcal{F}'$ which is an additive $\sigma$-skeleton, $X'_\alpha[n+1]$ is a posex of $X'_\alpha[n]$, so since $\mathfrak X$ satisfies the statement of our theorem, we can find an into isometry $g_\alpha[n+1]: X'_\alpha[n+1]\To \mathfrak X$ such that $g_\alpha[n+1]|_{X_\alpha[n]} = f^{-1}_\alpha[n]$. We define finally $X_\alpha[n+1] = g_\alpha[n+1](X_\alpha[n+1])$ and $f_\alpha[n+1] = g^{-1}_\alpha[n+1]$.\\

Proceeding this way, we have that $X_\alpha[n]\in\mathcal{F}$ for $n$ odd, while $X'_\alpha[n]\in \mathcal{F}'$ for $n$ even. At the end, $X_{\alpha+1} = \overline{\bigcup_{n<\omega}X_\alpha[2n+1]}\in\mathcal{F}$ and $X'_{\alpha+1} = \overline{\bigcup_{n<\omega}X'_\alpha[2n]}\in\mathcal{F}'$. $\qed$\\

\subsection{Universality property}
Let $\mathfrak X$ denote the space in Theorem~\ref{existsuniqueE}.

\begin{thm}
If $X$ is a tightly $\sigma$-filtered Banach space with $dens(X)\leq \frak c$, then $X$ is isometric to a subspace of the space $\mathfrak X$.
\end{thm}

Proof: Let $\{X_\alpha: \alpha \leq\kappa\}$ be subspaces that witness that $X$ is tightly $\sigma$-filtered. By the proof of Lemma~\ref{POidealE}, we can suppose that $\kappa$ is the cardinal $dens(X)$.  Then, inductively we can extend a given isometric embedding $X_\alpha\To \mathfrak X$ to $X_{\alpha+1}\To \mathfrak{X}$, by the properties of $\mathfrak X$.$\qed$\\

\section{Compact spaces}\label{sectionStone}

Along this section and the subsequent ones, $\mathfrak X$ will always denote the Banach space in Theorem~\ref{existsuniqueE} and $\mathfrak B$ the Boolean algebra in Theorem~\ref{existsuniqueB}.\\

\subsection{The compact space $\mathfrak K$}

\begin{defn}
Suppose that we have a commutative diagram of continuous surjections between compact spaces,

$$\begin{CD}
 K @>f>> L\\
 @VgVV @VvVV\\
 S @>u>>R.
\end{CD}$$
We say that this is a pull-back diagram if the following conditions hold:
\begin{enumerate}
\item For every $x,y\in K$, if $x\neq y$, then either $f(x)\neq f(y)$ or $g(x)\neq g(y)$.
\item If we are given $x\in S$ and $y\in L$ such that $u(x) = v(y)$, then there exists $z\in K$ such that $f(z)=y$ and $g(z)=x$.
\end{enumerate}
\end{defn}

Again, the notion of pull-back is more general in category theory, and in particular pull-back diagrams of continuous functions which are not surjective can be defined, but for our purposes we restrict to the case defined above.\\

If we are given two continuous surjections $u:S\To R$ and $v:L\To R$ between compact spaces, we can always construct a pull-back diagram as above making $K = \{(x,y)\in S\times L : u(x) = v(y)\}$ and taking $f$ and $g$ to be the coordinate projections. Moreover any other pull-back diagram 
$$\begin{CD}
 K' @>f'>> L\\
 @Vg'VV @VvVV\\
 S @>u>>R.
\end{CD}$$
is homeomorphic to the canonical one by a homeomorphism $h:K\To K'$ with $f'h=f$ and $g'h=g$.

\begin{prop}\label{stoneduality}
A diagram of embeddings of Boolean algebras
 $$\begin{CD} A @>>> B\\
 @AAA @AAA\\
 R @>>> S\end{CD}$$
is a push-out diagram if and only if and only if the diagram of compact spaces obtained by Stone duality,
$$\begin{CD}
 St(A) @<<< St(B)\\
 @VVV @VVV\\
 St(R) @<<< St(S),
\end{CD}$$
is a pull-back diagram.
\end{prop}

Proof: Left to the reader. $\qed$\\

\begin{prop}\label{dualballpullback}
Let $Y$ be a Banach space and $X,S,R$ subspaces of $Y$ such that $Y=\mathbf{PO}_R[S,X]$. Then the diagram obtained by duality between the dual balls endowed with the weak$^\ast$ topology,
$$\begin{CD}
 B_{Y^\ast} @>>> B_{X^\ast}\\
 @VVV @VVV\\
 B_{S^\ast} @>>>B_{R^\ast}
\end{CD}$$
is a pull-back diagram.
\end{prop}

Proof: We can suppose that $Y= (X\oplus_{\ell_1} S) / V$ where $V=\{(r,-r) : r\in R\}$, and then $B_{Y^\ast} \subset B_{(X\oplus_{\ell_1} S)^\ast} = B_{Y^\ast}\times B_{S^\ast}$ and it is precisely the set of pairs which agree on $R$.$\qed$\\

\begin{defn}
A continuous surjection between compact spaces $f:K\To L$ is called posex if there exists a pull-back diagram of continuous surjections $$\begin{CD}
 K @>f>> L\\
 @VVV @VVV\\
 S @>>>R.
\end{CD}$$
with $R$ and $S$ metrizable compact spaces.
\end{defn}

\begin{defn}
A compact space $K$ is called pull-back generated if there exists a family $\{K_\alpha : \alpha\leq\xi\}$ of compact spaces and continuous surjections $\{f_\alpha^\beta: K_{\beta}\To K_\alpha : \alpha\leq\beta\leq\xi\}$ such that
\begin{enumerate}
\item $K_0$ is a singleton and $K_\xi = K$,
\item $f_\alpha^\alpha$ is the identity map on $K_\alpha$,
\item $f_\alpha^\beta f_\beta^\gamma = f_\alpha^\gamma$ for $\alpha\leq\beta\leq\gamma\leq\xi$,
\item $f_\alpha^{\alpha+1}: K_{\alpha+1}\To K_\alpha$ is posex for every $\alpha<\xi$,
\item  If $\gamma\leq\xi$ is a limit ordinal and $x,y\in K_\gamma$ with $x\neq y$, then there exists $\beta<\gamma$ such that $f^\gamma_\beta(x)\neq f^\gamma_\beta(y)$ (this means that $K_\gamma$ is the inverse limit of the system below $\gamma$).
\end{enumerate}
\end{defn}

We can prove again a theorem about existence and uniqueness of a compact space in a similar way as we did for Banach spaces and Boolean algebras. But it is not worth to repeat the procedure because we would obtain just the Stone compact space of the Boolean algebra $\mathfrak B$ of Theorem~\ref{existsuniqueB}. Hence, we just denote this Stone space by $\mathfrak K = St(\mathfrak B)$.

\begin{lem}\label{zerodimensional}
Let $K$ be a pull-back generated compact space. Then there exists a zero-dimensional pull-back generated compact space $L$ of the same weight as $K$ and such that there is a continuous surjection from $L$ onto $K$.
\end{lem}

Proof: Assume that we have an inverse system $\{f_\alpha^\beta: K_\beta\To K_\alpha\}_{\alpha\leq\beta\leq\xi}$ as above witnessing that $K$ is pull-back generated. We produce our compact space $L$ and the continuous surjection $h:L\To K$ by constructing inductively a similar inverse system $\{g_\alpha^\beta:L_\beta\To L_\alpha\}_{\alpha\leq\beta\leq\xi}$ together with continuous surjections $h_\alpha:L_\alpha\To K_\alpha$ satisfying $h_\alpha g_\alpha^\beta = f_\alpha^\beta h_\beta$ for $\alpha\leq\beta$. We need that each $L_\alpha$ is zero-dimensional and the weight of $L_\alpha$ equals the weight of $K_\alpha$. The key step is providing $L_{\alpha+1}$, $f_\alpha^{\alpha+1}$ and $h_{\alpha+1}$ from $L_\alpha$ and $h_\alpha$. Since $f_\alpha^{\alpha+1}$ is posex there exist metrizable compact spaces $S$ and $R$ and a pull back diagram
$$\begin{CD}
 K_{\alpha+1} @>>> K_{\alpha}\\
 @VVV @VVV\\
 S @>>>R
\end{CD}$$
Consider $S'$ a metrizable zero-dimensional compact space and $u:S'\To S$ a continuous surjection. We have a larger diagram
$$\begin{CD}
@. @. L_\alpha\\
@. @. @V h_\alpha VV\\
 @. K_{\alpha+1} @>>> K_{\alpha}\\
 @. @VVV @VVV\\
S'@>u>> S @>>>R.
\end{CD}$$

We can define $L_{\alpha+1}$ and $g_\alpha^{\alpha+1}$ by making the pull back of the larger square above,
$$\begin{CD}
L_{\alpha+1} @>g^{\alpha+1}_\alpha>>  L_\alpha\\
@VVV @VVV\\
 S' @>>>  R.
\end{CD}$$

The continuous surjection $h_{\alpha+1}:L_{\alpha+1}\To K_{\alpha+1}$ can be obtained by applying the so called universal property of pull-back. In this case, the pull-back $K_{\alpha+1}$ can be seen as a subspace of $S\times K_\alpha$ and similarly $L_{\alpha+1}\subset S'\times L_{\alpha}$. One can define simply $h_{\alpha+1}(s,x) = (u(s),h_\alpha(x))$.$\qed$\\

\begin{prop}
For any diagram of continuous surjections between compact spaces
$$\begin{CD}
K\\
@VVV\\
 L @<<<  \mathfrak K,
\end{CD}$$
if $f:K\To L$ is posex and $weight(L)<\frak c$, then there exists a continuous surjection $\mathfrak K\To K$ that makes the diagram commutative.
\end{prop}

Proof: Without loss of generality, we suppose that we have a pull-back diagram
$$\begin{CD}
 K @>>> L\\
 @VVV @VVV\\
 S @>>>R
\end{CD}$$
where $S$ is metrizable and zero-dimensional. We can factorize into continuous surjections $\mathfrak K\To L'\To L$ such that $L'$ is zero-dimensional and $weight(L') = weight(L)$. Consider then $K'$ the pull-back of $K$, $L$ and $L'$, which is zero-dimensional because it is also the pull-back of $S$, $R$ and $L'$,
$$\begin{CD}
 S @<<< K @<<< K'\\
 @VVV @VVV @VVV \\
 R @<<< L @<<< L' @<<< \mathfrak K.
\end{CD}$$
 Then we have a similar diagram as in the statement of the theorem but all compact spaces are zero-dimensional. By Stone duality and property (3) of $\mathfrak B$ in Theorem~\ref{existsuniqueB} there is a continuous surjection $\mathfrak K\To K'$ that completes the diagram, and this provides the desired $\mathfrak K \To K$.$\qed$\\

\subsection{The relation between $\mathfrak X$ and $C(\mathfrak K)$}

\begin{thm}
$\mathfrak X$ is isometric to a subspace of $C(\mathfrak K)$.
\end{thm}

Proof: Let $\{X_\alpha : \alpha\leq\xi\}$ be an increasing chain of subspaces of $\mathfrak X$ witnessing the fact that $\mathfrak X$ is a push-generated Banach space. Then, by Proposition~\ref{dualballpullback} the dual balls $\{B_{X_\alpha^\ast} : \alpha\leq\xi\}$ form an inverse system that witness the fact that $B_{\frak X^\ast}$ is a pull-back generated compact space. By Lemma~\ref{zerodimensional} there exists a zero-dimensional pull-back generated compact space $L$ of weight $\mathfrak c$ that maps onto $B_{\mathfrak X^\ast}$. By Stone duality, using Proposition~\ref{stoneduality}, we get that the Boolean algebra $B$ of clopens of $L$ is a tightly $\sigma$-filtered Boolean algebra of cardinality $\mathfrak c$. Hence, by Proposition~\ref{Buniversal}, we can write $B\subset \mathfrak B$. Therefore $\mathfrak K$ maps continuously onto $L$, which maps continuously onto $B_{\mathfrak X^\ast}$. And this implies that $\mathfrak X\subset C(B_{\mathfrak X^\ast})\subset C(L) \subset C(\mathfrak K)$.$\qed$

\subsection{P-points}

Remember that a point $p$ of a topological space $K$ is called a $P$-point if the intersection of countably many neighborhoods of $p$ contains a neighborhood of $p$. The following result was proven by Rudin~\cite{Rudin} under CH and by Stepr\={a}ns~\cite{Steprans} in the $\aleph_2$-Cohen model: for every two $P$-points $p,q\in\omega^\ast$ there exists a homeomorphism $f:\omega^\ast\To \omega^\ast$ such that $f(p)=q$. Geschke~\cite{Geschke} proves that it is sufficient to assume that $\mathcal{P}(\omega)$ is tightly $\sigma$-filtered. In this section we prove that this is a property of our compact space $\mathfrak K$.

\begin{thm}\label{ppoints}
Let $p,q\in\mathfrak K$ be P-points. Then there exists a homeomorphism $f:\mathfrak K \To \mathfrak K$ such that $f(p) = q$.
\end{thm}

In what follows, points of the Stone space of a Boolean algebra $B$ are considered as ultrafilters on $B$. Given $Q_1,Q_2\subset B$ we write $Q_1\leq Q_2$ if $q_1\leq q_2$ for every $q_1 \in Q_1$ and $q_2 \in Q_2$. If $b\in B$ and $Q\subset B$, then $b\leq Q$ means $\{b\}\leq Q$ and $Q\leq b$ means $Q\leq \{b\}$.

\begin{defn}If $A$ is a subalgebra of $B$, $Q\subset A$ and $b\in B$, then we write $b\leq_A Q$ if $\{a\in A : b\leq a\} $ equals the filter generated by $Q$ in $A$. Similarly, we write $Q\leq_A b$ if $\{a\in A : a\leq b\}$ equals the ideal generated by $Q$ in $A$. \end{defn}

With this notation, the notion of push-out of Boolean algebras can be rephrased as follows: given Boolean algebras $R\subset A,S\subset B$, then $B=\mathbf{PO}_R[S,A]$ if and only if for every $b\in S$, $\{a\in R : a\leq b\}\leq_A b \leq_A \{a\in R : b\leq a\}$.

\begin{lem}\label{ppointimproved}
Let $p\in\mathfrak K$ be a P-point, let $A\subset \mathfrak B$ be a subalgebra with $|A|<\frak c$ and let $Q\subset A\cap p$ be countable. Then, there exists $b\in p$ such that $\{0\} \leq_A b \leq_A Q$.
\end{lem}

Proof: Because $p$ is a $P$-point and $Q$ is countable, there exists $b_0\in p$ such that $b_0\leq Q$. Since $|A|<\mathfrak c$, any posex extension of $A$ can be represented inside $\mathfrak B$. In particular, we can find $b_1\in\mathfrak B$ such that $\{0\}\leq_A b_1\leq_A Q$  in such a way that $B_1 = \langle A\cup\{b_1\}\rangle$ is a posex of $A$. Observe that $b_0\cup b_1\leq_A Q$ and $b_0\cup b_1\in p$. By the same reason as before, we can find $b_2\in \mathfrak B$ such that $\{0\}\leq_{B_1} b_2 \leq_{B_1} \{1\}$.  Notice that we have $\{0\}\leq_A (b_1\cup b_0)\cap b_2 \leq_A Q$ and $\{0\}\leq_A (b_1\cup b_0)\setminus b_2\leq_A Q$ and since $b_1\cup b_0\in p$ and $p$ is an ultrafilter, either $(b_1\cup b_0)\cap b_2$ or $(b_1\cup b_0)\setminus b_2$ belong to $p$.$\qed$\\

\begin{lem}\label{ppointextension}
Let $A$ be a Boolean algebra with $|A|<\frak c$, and $u:A\To \mathfrak B$ and $v:A\To B$ be Boolean embeddings with $v$ being posex. Fix $P$-points, $p_0$, $p$ and $q$ of $A$, $\mathfrak B$ and $B$ respectively such that $u(p_0)\subset p$ and $v(p_0)\subset q$. Then there exists an embedding $\tilde{u}:B\To \mathfrak B$ such that $\tilde{u}v = u$ and $\tilde{u}(q)\subset p$.
\end{lem}

Proof: Suppose that $B = \mathbf{PO}_R[S,v(A)]$ where $S$ is a countable subalgebra of $B$. We can produce a further posex superalgebra $B_0\supset B$ of the form $B_0 = \langle B\cup\{b_0\}\rangle$ such that $\{0\}\leq_B b_0 \leq_B S\cap q$. By Lemma~\ref{ppointimproved} we can find $b_1\in p$ such that $\{0\}\leq_{u(A)} b_1\leq_{u(A)} u(R\cap p_0)$. Notice that $\{0\} \leq _{v(A)} b_0 \leq_{v(A)} v(R\cap p_0)$, and this allows to define a Boolean embedding $w:\langle v(A)\cup\{b_0\}\rangle \To \mathfrak B$ such that $w v = u$ and $w(b_0) = b_1$. Since $B_0$ is a posex of $v(A)$ and $\langle v(A)\cup\{b_0\}\rangle$ is countably generated over $v(A)$, we have that $B_0$ is a posex of $\langle v(A)\cup\{b_0\}\rangle$, so using the second property stated in Theorem~\ref{existsuniqueB}, we find $\tilde{w}:B_0\To \mathfrak B$ such that $\tilde{w}|_{\langle v(A)\cup\{b_0\}\rangle} = w$. We consider finally $\tilde{u} = \tilde{w}|_B$. It is clear that $\tilde{u}v = u$. On the other hand, since $\tilde{u}(b_0) = w(b_0) = b_1\in p$ and $b_0\leq S\cap q$, we have that $\tilde{u}(S\cap q) \subset p$. It is also clear that $\tilde{u}(q\cap v(A)) = u(p_0) \subset p$. Since $B = \langle S\cup v(A)\rangle$ the ultrafilter $q$ is the filter generated by $(q\cap S)\cup(q\cap v(A))$, so we finally get that $\tilde{u}(q)\subset p$.$\qed$\\

We can rephrase the statement of Theorem~\ref{ppoints} as follows: ``If $\mathfrak B$ and $\mathfrak B'$ are Boolean algebras satisfying the conditions of Theorem~\ref{existsuniqueB} and $p$ and $q$ are $P$-points in $\mathfrak B$ and $\mathfrak B'$ respectively, then there exists an isomorphism $f:\mathfrak B\To \mathfrak B'$ such that $f(p) = q$''. In order to prove this, we just have to follow the proof of uniqueness in Theorem~\ref{existsuniqueB} and make sure that at each step it is possible to choose the partial isomorphisms $f_\alpha: B_\alpha\To B'_\alpha$ in such a way that $f(p\cap B_\alpha) = q\cap B'_\alpha$. And what we need for that is exactly Lemma~\ref{ppointextension} applied to $\mathfrak B$ and $\mathfrak B'$.

\section{Open problems}\label{sectionProblems}

\subsection{When $\mathfrak c$ is singular}

\begin{prob}
Do Theorem~\ref{existsuniqueE} and Theorem~\ref{existsuniqueB} hold when $\mathfrak c$ is singular? 
\end{prob} 

Here the point is that regularity looks essential to control all substructures of size less than $\mathfrak c$ in the \emph{existence} part of the proof. Perhaps these theorems do not hold for singular $\mathfrak c$ as they are stated, but it would be satisfactory any variation that would allow us to speak about unique objects $\mathfrak X$ and $\mathfrak B$ defined by some properties in ZFC. We remark that if we restrict below a given cardinal, \emph{existence} can be proven. We state it for Boolean algebras, and leave the Banach space version to the reader.

\begin{prop}
Fix a cardinal $\lambda<\mathfrak c$. There exists a Boolean algebra $\mathfrak B$ with the following properties:
\begin{enumerate}
\item $|\mathfrak B| = \mathfrak c$,
\item $\mathfrak B$ is tightly $\sigma$-filtered,
\item For any diagram of embeddings of the form
$$\begin{CD}
 B \\
 @AAA \\
 A @>>> \mathfrak B,
\end{CD}$$
if $|A|\leq\lambda$ and $A\To B$ is posex, then there exists an isometric embedding $B\To \mathfrak B$ which makes the diagram commutative.
\end{enumerate}
\end{prop}

Proof: Let $\lambda^+$ be the successor cardinal of $\kappa$. Construct $\mathfrak B$ in the same way as in the proof of Theorem~\ref{existsuniqueB} or Theorem~\ref{existsuniqueE}, but do it with a tower of length $\mathfrak c \cdot \lambda^{+}$ instead of length $\mathfrak c$.$\qed$\\

\subsection{Properties of $\mathfrak B$ reflected on Banach spaces}

We know much more about $\mathfrak B$ than about its Banach space relative $\mathfrak X$. So it is natural to ask whether certain facts that hold for $\mathfrak B$ in some models are reflected by analogous properties for $\mathfrak X$ or for $C(\mathfrak K)$.\\

For example, we mentioned that in the $\aleph_2$-Cohen model, $\mathfrak B = \mathcal{P}(\omega)/fin$ \cite{DowHart}. In this model, therefore, $\mathfrak B$ has some additional properties of extensions of morphisms: given any diagram of arbitrary morphisms
$$\begin{CD}
 S\\
 @AAA \\
 R @>>>\mathfrak B
\end{CD}$$
where $R$ is countable and $S$ is arbitrary, there exists a morphism $S\To\mathfrak B$ which makes the diagram commutative. The analogous property for Banach spaces is called (1-)universally separably injectivity (we add a 1 if the operator $S\To\mathfrak X$ can be found with the same norm as $R\To\mathfrak X$), cf.~\cite{extremadura}. This property implies that the space contains $\ell_\infty$.

\begin{prob}
In the $\aleph_2$-Cohen model, is the space $\mathfrak X$ universally separably injective? Does it contain $\ell_\infty$? 
\end{prob}

Another observation is that in this model, $\mathcal{P}(\omega)/fin$ does not contain any $\omega_2$-chain. Indeed, Dow and Hart~\cite{DowHart} prove in ZFC that no Cohen-Parovi\v{c}enko Boolean algebra (in particular $\mathfrak B$) can contain $\omega_2$-chains. On the other hand, Brech and Koszmider~\cite{BrechKoszmider} prove that in the $\aleph_2$-Cohen model, $\ell_\infty/c_0$ does not contain the space $C[0,\omega_2]$ of continuous functions on the ordinal interval $[0,\omega_2]$. So the natural question is:

\begin{prob} Is it true in general that $C(\mathfrak K)$ (or at least $\mathfrak X$) does not contain $C[0,\omega_2]$?
\end{prob}

\begin{prob}
Can spaces like $\ell_2(\omega_2)$ or $c_0(\omega_2)$ be subspaces of $\mathfrak X$? 
\end{prob}

Under (MA + $\mathfrak c = \aleph_2$), Dow and Hart \cite{DowHart} prove that $\mathfrak B$ does not contain $\mathcal{P}(\omega)$ as a subalgebra. So we may formulate

\begin{prob}[MA + $\mathfrak c = \aleph_2$] Does $C(\mathfrak K)$ contain a copy of $\ell_\infty$?
\end{prob}

This question is already posed in \cite{extremadura}. A negative answer would solve a problem by Rosenthal by providing an $F$-space $K$ such that $C(K)$ does not contain $\ell_\infty$. In~\cite{extremadura} it is proven that, under MA + $\mathfrak c = \aleph_2$, there is an isometric embedding $c_0\To C(\mathfrak K)$ which cannot be extended to an embedding $\ell_\infty\To C(\mathfrak K)$.\\

The space $\mathfrak X$ plays the role of $\mathfrak B$ in the category of Banach spaces. But we do not have a Banach space playing the role of $\mathcal{P}(\omega)/fin$, when it is different from $\mathfrak B$. This is related to the general question of finding intrinsic characterizations of $\mathcal{P}(\omega)/fin$ out of CH and $\aleph_2$-Cohen models, that we could translate into other categories.  

\begin{prob}
Is there a Banach space counterpart of $\mathcal{P}(\omega)/fin$?
\end{prob}

\section{Other cardinals}\label{sectionCardinals}

We comment that a more general version of some of our results can be stated if we let arbitrary cardinals to play the role of $\mathfrak c$ and countability. We state it for Banach spaces, and leave the Boolean version to the reader. For an uncountable cardinal $\tau$, we say that an isometric embedding $X\To Y$ is $\tau$-posex if $Y = \mathbf{PO}[S,X]$ for some $S$ with $dens(S)<\tau$. We say that $X$ is tightly $\tau$-filtered if $X$ is the union of a continuous tower of subspaces starting at 0 and such that each $X_{\alpha+1}$ is a $\tau$-posex of $X_\alpha$ (in the case of Boolean algebras, such a definition is equivalent to the one given in \cite{Geschke}).

\begin{thm}\label{existsuniqueEcardinal} Let $\kappa$ be a regular cardinal, and $\tau$ an uncountable cardinal with $\kappa^{<\tau} = \kappa$. 
There exists a unique (up to isometry) Banach space $\mathfrak X = \mathfrak{X}(\kappa,\tau)$ with the following properties:
\begin{enumerate}
\item $dens(\mathfrak X) = \kappa$,
\item $\mathfrak X$ is tightly $\tau$-filtered,
\item For any diagram of isometric embeddings of the form
$$\begin{CD}
 Y \\
 @AAA \\
 X @>>> \mathfrak X,
\end{CD}$$
if $dens(X)<\kappa$ and $X\To Y$ is $\tau$-posex, then there exists an isometric embedding $w:Y\To \mathfrak X$ which completes the diagram.
\end{enumerate}
\end{thm}

The proof would be just the same. Our original space corresponds to $\mathfrak X(\mathfrak c,\aleph_1)$.\\

\end{document}